\documentclass{article}
\usepackage[utf8]{inputenc}
\usepackage{cite}
\usepackage{amsmath,amssymb,amsfonts}
\usepackage{mathabx}
\usepackage{algorithmic}
\usepackage{graphicx}
\usepackage{textcomp}
\usepackage{xcolor}

\usepackage{hyperref}

\usepackage{amsthm}
\theoremstyle{plain}

\providecommand{\keywords}[1]
{
  \small	
  \textbf{\textit{Keywords---}} #1
}

\title{Graph Polynomials  for the Numbers of Independent Sets and Bipartite Cuts for Undirected Graphs}
\author{R. L. Streit$^{1}$ \\
\small $^1$Metron, 1818 Library St., Suite 600,  Reston, VA 20190  USA\\
\\
\small 
Corresponding author. E-mail: streit@metsci.com} 

\date{}

\begin{document}

\maketitle

%\tableofcontents

\abstract {%Graph polynomials  are derived for  four   problems in the theory of general undirected graphs.  %, but they are derived as     coefficients of a multivariate GF of the number of solutions of appropriately  defined systems of linear  Diophantine equations.  
The  graph polynomial for the number of independent sets of size $k$ in a general undirected graph is  shown to be equal to an elementary symmetric polynomial of the vertex monomials, which are determined by the edges incident at the vertices.   The edge variables that comprise the vertex monomials are shown to be nilpotent of degree two.  The index of nilpotency of the algebra generated by the graph's vertex monomials is the  order  of a maximum independent set of the graph. 
Graph polynomials for the number of cliques and vertex covers are derived  from the graph polynomial for independent sets. %   using  well-known relationships between these concepts.  
The graph polynomial for the number of bipartite cuts is given as a coefficient of  a multivariate Laurent polynomial.  

\bigskip

\keywords{Independent sets; graph polynomials; %graph invariants; 
nilpotent algebra; nilpotent index; %cliques; vertex covers; 
bipartite partitions; Laurent polynomials}

\section{Introduction}
\label{Introduction}
Graph polynomials are univariate  generating functions (GFs) for graph invariants.  Examples of  graph polynomials for general undirected graphs include  independence polynomials  \cite{levitmandrescu, GPoly1}, chromatic polynomials \cite[Ch. 9]{biggs2ed}, Tutte polynomials \cite[Ch. 13]{biggs2ed}, matching polynomials \cite{GPoly2}, flow polynomials \cite{GodsilRoyle}, and circuit partition polynomials \cite{EllisMonaghan} \cite{bollaCPP}.   
Graph polynomials  for  the  numbers of independent sets, cliques, vertex covers, and bipartite cuts are derived in this paper. The independence polynomial is derived first.

The independence polynomial is the graph polynomial for the  number of independent sets.   We formulate it as a problem of counting the number of nonnegative integer solutions of a system of linear equations.   Lattice point counting methods (see, e.g., \cite{deloera}) can be used, but we  take a different approach, one that encodes the counting problem into a multivariate analytic function. We call this function a partition function because it is defined as a power series summed over the set of solutions of the system (\textit{cf.}  Eqn. \eqref{multiGF} below). As is the case in many successful applications of analytic combinatorics, so it is here also---the partition function for graph independence polynomials takes a simple analytic form that can be written down directly from the graph structure.  The  independence polynomial is, thus,  the coefficient of a specific monomial in the (multivariate) power series expansion of the partition function.    

Using this approach we show that the independence polynomial  is  an elementary symmetric polynomial (ESP) in the vertex monomials of the graph (\textit{cf.}  Eqn. \eqref{vertexmonomial} below). The derivation shows that the edge variables that comprise the vertex monomials  are nilpotent of degree two, that is, the square of each variable is zero. Nilpotency significantly reduces the computational complexity of  evaluating the ESPs.  The algebra of the graph vertex monomials is nilpotent, and the index of nilpotency is equal to one plus the size of the maximum independent set. Computing the maximum independent set is an  NP-complete problem, so  it follows that computing the nilpotent index of an algebra is also NP-complete.  

The  graph polynomials for  cliques and vertex covers are derived from the independence polynomial. These polynomials  give  a different perspective on the relationships between these concepts.  The graph polynomial for bipartite cuts is also derived using partition functions, but in this case the  partition function is a Laurent series.

\textit{Scholium.} It is stated in \cite{Rivin} that  Euler \cite{Euler} was the first to solve systems of linear Diophantine equations using the  method that we use here.  For further discussion, \cite{Rivin} refers readers to   \cite[p. 375-6]{Schrijver}.

\section{Background and outline}
A general undirected  graph is denoted by  $\mathbf G=(V,E)$, where $V$ is the set of vertices and $E$ is the set of edges.  The number of vertices, $n\equiv \left| V\right|$, is called the order of the graph, while the number of edges, $m\equiv \left|E\right|$, is called the size of the graph.  Throughout this paper,  $V$ and $E$ are assumed to be finite, and  $\mathbf G$ is assumed to be simple, meaning that there is at most one edge connecting any two vertices, and that there are no loops (i.e.,  no edge  connects a vertex to itself). 

An independent set of vertices of the graph $\mathbf G$ is a subset of $V$  such that no two  vertices in the subset are adjacent, i.e., they are not connected by an edge in $E$.  The order  of an independent set is its cardinality. An independent set is  maximal if it is not a subset of any other  independent set. A maximum independent set is an independent set of the largest possible order. Every maximum independent set is maximal, but not conversely; e.g., the  3-cube graph has six   maximal  sets---four of order two and two of order four, the maximum order.  

Let  $A_k$ denote the number of independent sets of  order $k\ge 0$.  The GF of the sequence $A_0, A_1, A_2, \ldots$ defines the graph independence polynomial, 
\begin{align}
\label{GFis}
\pi_{\alpha}(z)=\sum_{k=0}^\infty A_k z^k,
\end{align}
where $z$ is a real or complex valued  indeterminate variable. Because  $V$ is finite,  $\pi_{\alpha}(z)$ is a polynomial whose degree, denoted by $\alpha(\mathbf G)$,  is the independence number, i.e., the order of a maximum independent set. %The GF  \eqref{GFis} is  the graph polynomial for the number of independent sets.   

Section \ref{examples} introduces the partition function,  $\Psi_{\mathbf G}(\cdot)$, that we use to derive $\pi_{\alpha}(z)$.  We begin the section  with a simple example. The technique extends to the clique, vertex cover, and bipartite cut problems. 

Section \ref{explicitISform} derives  an explicit expression for  $\pi_{\alpha}(z)$ for  general undirected graphs.  The coefficients $A_k$ are shown to be  elementary symmetric polynomials in the vertex monomials, where a vertex  monomial is defined by the edges  incident to the vertex (\textit{cf.}  Eqn. \eqref{vertexmonomial}). 
Stated differently, the result links  the independent set problem  to a nilpotent algebra of the vertex monomials of the graph.      

A clique is a subset of  $V$ such that every pair of distinct vertices in the subset are adjacent.  The graph polynomial of the number of cliques in $\mathbf G$ is derived in Section \ref{Scetioncliques} from the independence polynomial using the fact that the complement of a clique is an independent set, and conversely.     
    
Section \ref{vertexcovers} defines a vertex cover as a subset of  $V$ that includes at least one endpoint of every edge in $E$. The  order of a cover is its cardinality. The graph polynomial  of the number of vertex covers is derived from the  polynomial  of the number of  independent sets.   

Section \ref{maxcut} defines bipartite cut-sets. A cut-set  comprises a set of edges that, when removed,  partition $V$ into two disjoint parts.  The size of a cut-set is the number of edges. The problem is modeled using  the XOR (exclusive-or).  The graph polynomial  of the number of  cut-sets  of a given size is derived from an appropriately defined partition function.   

Section \ref{Concludingremarks} briefly comments on using  ``labeled'' partition functions  to enumerate  a  complete list of the independent sets, cliques, vertex covers, and bipartite cuts.  Weighted versions of these  problems, e.g.,   Ising problems from physics, are also mentioned. The section concludes with a brief discussion of asymptotic  approximations.   
The Hardy-Littlewood circle method and, more generally  saddle point methods \cite{FS}, are  applicable.

\section{Independence polynomials }
\label{examples}

.  

The maximum independent set problem for  undirected graphs is often posed (e.g., \cite{princetonUni}) as a  binary linear programming (BLP) problem.  We  modify the BLP formulation  to find the coefficients of  $\pi_{\alpha}(z)$.  

\newtheorem{lemma}{Lemma}
\begin{lemma}
For each order $k\ge 0$,  define  the following binary (Diophantine) system of $m+1$ equations in $n+m$ variables:   
\begin{subequations} 
\label{MISSlackForm}
\begin{align}
&\textstyle \sum_{v\in V} x_v = k\label{MISSlackForm2a}\\
&x_u+x_v +x_{uv} = 1,\,\, \textrm{for all }  uv\in E, \label{MISSlackForm2b}\\
\textrm{where}\quad &x_v \in \{0,1\}\,\, \textrm{and}\,\,
x_{uv}\in \{0,1\}.\label{MISSlackForm2d}
\end{align}
\end{subequations}
The number of $\{0,1\}$ solutions of this system is $A_k$, the number of independent sets of  order $k\ge 0$. 
\end{lemma}
To see this we show that  there is a one-to-one correspondence between the set of solutions of \eqref{MISSlackForm} and the  collection of independent sets of the graph $\mathbf G$ of order $k$. Given   $V^*$, an independent set of order $m$, define  $x_v^*=1$ if $v\in V^*$, and $x_v^*=0$ if not.  Let $\mathbf x^*_V=\{x_v^* : v\in V\}$.  Given $\mathbf x^*_V$, the edge variables $\mathbf x_E^*=\{x_{uv}^* : uv\in E\}$ are uniquely defined by \eqref{MISSlackForm2b}.  Thus, $V^*$ maps to a  solution of \eqref{MISSlackForm}.  Moreover,  no two independent sets of order $k$ map to the same solution.   Conversely,  given a solution  $\mathbf x_V^*$ and $\mathbf  x_{E}^*$  of  \eqref{MISSlackForm}, let $V^* = \{v\in V : x_v^* =1\}$.   From \eqref{MISSlackForm2a}, the order of $V^*$ is $k$. From \eqref{MISSlackForm2b}, the vertices in $V^*$ are not adjacent, so $V^*$ is an independent set. Given $\mathbf x_V^*$,  the  edge variables $\mathbf x_E^*$ are uniquely defined, so different solutions of  \eqref{MISSlackForm} map to different independent sets of order $k$. 
\qed

\bigskip
We use this result to count the number of solutions of the Diophantine system. This derivation gives us an expression for the independence polynomial as the coefficient of a multivariate polynomial. As will be shown, the product form of this polynomial can be written   easily from the edge-to-vertex incidence matrix of the graph.  

\subsection{Example: counting the number of solutions}
\label{Example1}
The graph depicted in Fig. \ref{sixnodeexample} comprises  
\begin{figure}[b]
\centering 
\includegraphics[width=0.5\textwidth]{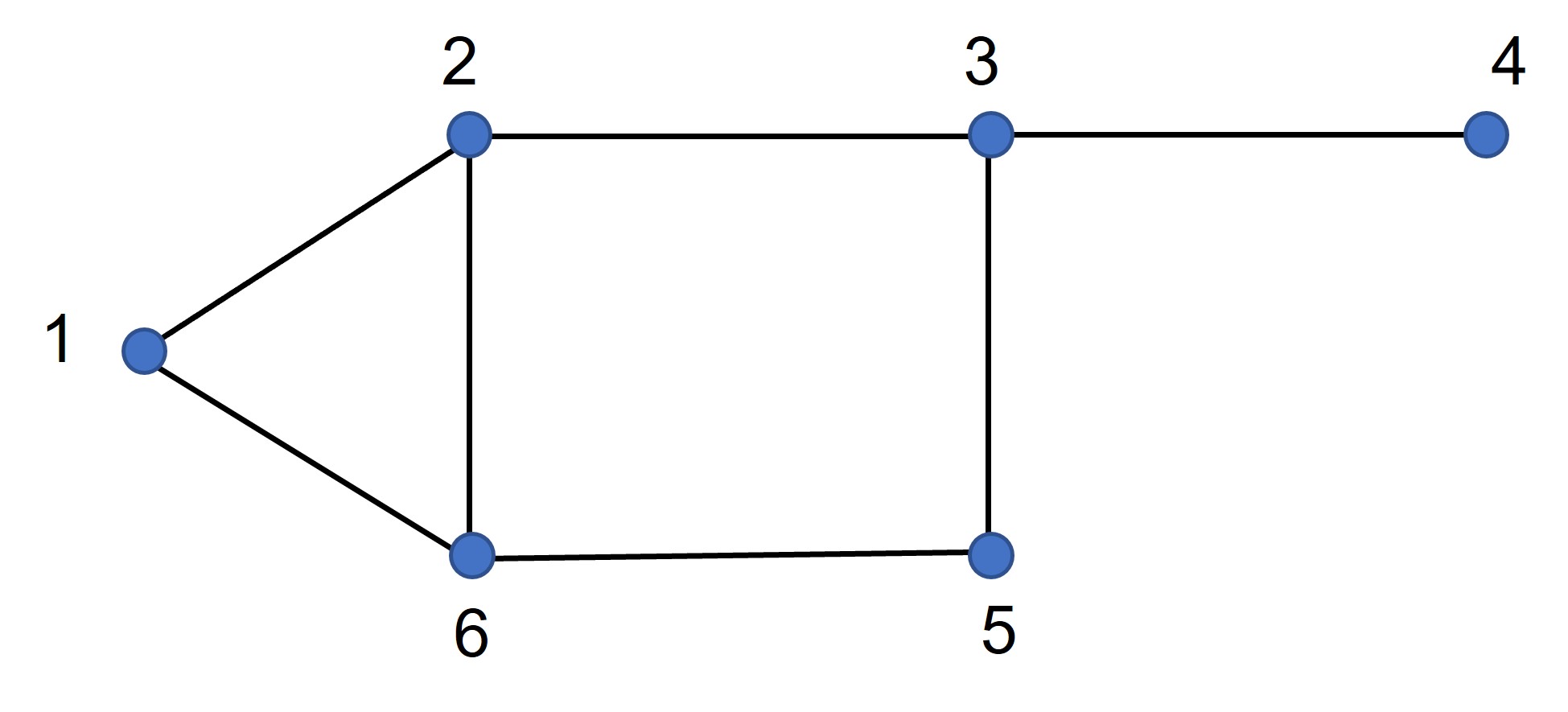}
    \caption{Undirected graph example}
\label{sixnodeexample}
\end{figure}  
a total of 13 variables, one for each of the  $n=6$  vertices and  $m=7$ edges.   Let   
\begin{align}
%\label{orderly}
 \mathbf x&=(\mathbf x_V,\mathbf  x_E)\nonumber\\
 &\equiv \big(\underbrace{x_1,x_2,x_3,x_4,x_5,x_6}_6,\underbrace{x_{12},x_{16},x_{23},x_{26},x_{34},x_{35},x_{56}}_7\big)\in\{0,1\}^{13}.\nonumber
\end{align}
The first six  are the ``vertex'' variables $\mathbf x_V$, and the others are the ``edge'' variables $\mathbf x_E$.  The vector $\mathbf x$ imposes a linear ordering on the variables assigned to the vertices and edges.  The  equations  \eqref{MISSlackForm2a}-\eqref{MISSlackForm2b} are 
\begin{equation}
\label{exampleequ}
\begin{split}
x_{1}+x_{2}+x_{3}+x_{4}+x_{5}+x_{6}&=k\\
    x_1+x_2+x_{12}&=1\\
    x_1+x_6+x_{16}&=1\\
    x_2+x_3+x_{23}&=1\\
    x_2+x_6+x_{26}&=1\\
    x_3+x_4+x_{34}&=1\\
    x_3+x_5+x_{35}&=1\\
    x_5+x_6+x_{56}&=1.
\end{split}
\end{equation}
Assign the indeterminant   $z$ to the first equation  and the indeterminants $z_{12},\ldots,z_{56}$ to the other equations. Indeterminate variables are, here, independent real or complex variables.   The vector $\mathbf x$ is a solution of the system if and only if   
\begin{align}
   &z^{ x_{1}+x_{2}+x_{3}+x_{4}+x_{5}+x_{6}}      \,z_{12}^{x_1+x_2+x_{12}} \,z_{16}^{x_1+x_6+x_{16}}\cdots \,z_{56}^{x_5+x_6+x_{56}} \nonumber\\
    &\qquad=z^k z_{12} z_{16}\cdots z_{56}\, .\nonumber  
\end{align}
Rearranging  terms gives the equivalent equation 
\begin{align}
\label{examlabel}
\begin{split}    &(z\,z_{12} z_{16})^{x_{1}}(z \,z_{12} z_{23} z_{26})^{x_{2}}\cdots (z\,z_{16} z_{26} z_{56})^{x_{6}}\,z_{12}^{x_{12}}\,z_{16}^{x_{16}}\cdots z_{56}^{x_{56}}\\
    &\qquad =z^k z_{12} z_{16}\,\cdots\,z_{56}. 
\end{split}
\end{align}
The number of nonnegative integer solutions $\mathbf x$ of \eqref{examlabel} is the number of nonnegative integer solutions of \eqref{exampleequ}, which is the number of independent sets of order  $k$. Let  $\mathbf z_E=(z_{12},\ldots,z_{56})$, and define the polynomial 
\begin{align}
    \Psi_{\mathbf G}(z, \mathbf z_E)=\sum_{\mathbf x\in\{0,1\}^{13}}
\bigg\{\begin{split}
   & (z\,z_{12} z_{16})^{x_{1}}(z \,z_{12} z_{23} z_{26})^{x_{2}} \\
    & \qquad\cdots\; (z\,z_{16} z_{26} z_{56})^{x_{6}}\,z_{12}^{x_{12}}\,z_{16}^{x_{16}}\cdots z_{56}^{x_{56}}
    \end{split}\bigg\}. \nonumber
\end{align}
Note that this polynomial is independent of $k$. The sum is separable, so it factors into  a product of 13 terms:  
\begin{equation}
\label{neededattheend}
   \begin{split} \Psi_{\mathbf G}(z, \mathbf z_E)&=(1+z\, z_{12} z_{16}) (1+z\, z_{12} z_{23} z_{26})\cdots(1+z\, z_{16} z_{26} z_{56}) \\
    &\qquad \times (1+z_{12})(1+z_{16})\cdots(1+z_{56}).
    \end{split}
\end{equation} 
Expanding this product into a multivariate power series about the origin in $\mathbb C^{m+1}$ gives  
\begin{align*}
    \Psi_{\mathbf G}(z, \mathbf z_E)&=\sum_{k_0,k_{12},\ldots,k_{56}=0}^\infty a(k_0,k_{12},k_{16},\ldots,k_{56})\;z^{k_0}z_{12}^{k_{12}}z_{16}^{k_{16}}\cdots z_{56}^{k_{56}},
\end{align*}
where the coefficients $a(k_0,k_{12},k_{16},\ldots,k_{56})$ are nonnegative integers.   It is now self evident that the coefficient of the monomial  $z^k z_{12}z_{16}\cdots z_{56}$ is equal to the number of solutions of \eqref{exampleequ}, that is, 
\[
A_k=a(k,1,1,\ldots,1).
\] 
This result is more clearly written in  bracket notation, which  is widely used in analytic combinatorics to denote specific   coefficients of  (possibly formal) power series expansions.    
For example, $c_{jk} = [x^j y^k]f(x,y)$ denotes the coefficient of $x^jy^k$ in the  series   $f(x,y)=\sum_{j,k\ge 0} c_{jk} x^j y^k$.  The notation extends  to any number of variables in the obvious manner.  

Using  bracket notation gives  
\begin{align}
    A_k &=\big [z^k z_{12}z_{16}\cdots z_{56} \big ] \Psi_{\mathbf G}(z, \mathbf z_E),\quad k\ge 0\nonumber.
\end{align} 
It follows that the  coefficient of $z_{12}z_{16}\cdots z_{56}$ is a polynomial in $z$, and that this polynomial is the GF of the number of independent sets:    
\begin{align}
\label{mypiExam}
    \pi_{\alpha}(z) &=\big [z_{12}z_{16}\cdots z_{56} \big ] \Psi_{\mathbf G}(z, \mathbf z_E)\nonumber\\
    &=1 + 6 z + 8 z^2 + 2 z^3.
\end{align}
It is seen from this GF  that there are $2$  maximum independent sets of order $\alpha(\mathbf G)=3$. By inspection, the maximum sets are $\{1,4,5\}$ and $\{2,4,5\}$. The next section treats the general case.

\subsection{General case }
\label{multivariate}
Let $\mathbf x_V\in \{0,1\}^{{n}}$ denote the vector of  vertex variables $x_v,\,v\in V$, arranged in some arbitrary, but fixed, linear order.  Similarly, we let $\mathbf x_E\in \{0,1\}^{{m}}$ denote the vector of edge variables $x_{uv},\,uv\in E$, in some fixed order.  
We assign the indeterminant  $z\in\mathbb C$ to  equation \eqref{MISSlackForm2a} and the vector of indeterminants $\mathbf z_E=\{z_{uv},\,uv\in E\}\in\mathbb C^m$ to the  equations  \eqref{MISSlackForm2b}.  The  components of  $\mathbf z_E$ are arranged in the same order  as $\mathbf x_E$.    The vector $(\mathbf x_V, \,\mathbf x_E)$ is a solution of  \eqref{MISSlackForm2a}-\eqref{MISSlackForm2b} if and only if  
\begin{align}
\label{EulerRocks}
   &z^{\sum_{v\in V} x_v}      \prod_{uv\in E} z_{uv}^{x_u+x_v+x_{uv}}=z^k\! \prod_{uv\in E} z_{uv}  
\end{align}
is an algebraic identity in the indeterminate variables $(z,\mathbf z_E)$. Note that $k$ appears only on the right hand side of the equation.   

We define the partition function of the graph $\mathbf G$ to be  the sum over all feasible  $\mathbf x_V$ and $\mathbf x_E$ of the left hand side: 
\begin{align}
    &\Psi_{\mathbf G}(z,\mathbf z_E)=\sum_{(\mathbf x_V\!, \mathbf x_E)\in\{0,1\}^{n+m}} \left( z^{\sum_{v\in V} x_v}      \prod_{uv\in E} z_{uv}^{x_u+x_v+x_{uv}}\right )\nonumber\\
    &\qquad \quad=\sum_{(\mathbf x_V\!, \mathbf x_E)\in\{0,1\}^{n+m}} \left( \prod_{v\in V} z^{x_v} \right) \!\left (     \prod_{uv\in E} z_{uv}^{x_u+x_v} \right) \!\left (\prod_{uv\in E} z_{uv}^{x_{uv}}\right ).\label{theGfatlast}
\end{align}
Let $E(v)\subset E$ denote the set of  edges incident to the vertex $v\in V$, and define the ``vertex monomials''  as  products of edge variables,  
\begin{align}
\label{vertexmonomial}
    Z_{E(v)}=\prod_{uv\in  E(v)} z_{uv},\;\text{ for } v\in V.
\end{align}
With this notation, rearranging the terms in \eqref{theGfatlast} gives 
\begin{align}
\label{neversayneverlabel}
    \Psi_{\mathbf G}(z,\mathbf z_E)&=\sum_{(\mathbf x_V\!, \mathbf x_E)} \left( \prod_{v\in V} z^{x_v}    \prod_{v\in V} Z_{E(v)}^{x_{v}}\prod_{uv\in E} z_{uv}^{x_{uv}}\right )\nonumber\\
    &=\sum_{(\mathbf x_V\!, \mathbf x_E)} \left( \prod_{v\in V} \left (z  \,  Z_{E(v)}\right )^{x_v}\!\prod_{uv\in E} z_{uv}^{x_{uv}}\right ).
\end{align}
The sum is separable in the  variables  ${\mathbf x_V}$ and ${\mathbf x_E}$.  Consequently,     
\begin{align}
\label{multiGF}
    \Psi_{\mathbf G}(z,\mathbf z_E)&=\prod_{v\in V}\left(\sum_{\mathbf x_v\in\{0,1\}} \left (z  \,  Z_{E(v)}\right )^{x_v}\right)\prod_{uv\in E}\left(\sum_{\mathbf x_{uv}\in\{0,1\}}z_{uv}^{x_{uv}}\right)\nonumber\\
    &=\prod_{v\in V}\left(1+z  \,  Z_{E(v)}\right)\!\prod_{uv\in E}\left(1+z_{uv}\right).
\end{align}
The graph is finite, so the partition function is a multivariate polynomial. 

\begin{lemma}
\label{lemma2}
    The independence polynomial of the graph $\mathbf G=(V,E)$ is the coefficient  \begin{align}
\label{mypie}
    \pi_{\alpha}(z) =\Big [\textstyle \prod_{uv\in E} z_{uv} \Big ] \Psi_{\mathbf G}\!\left(z,\mathbf z_E\right),
\end{align}
where  $\Psi_{\mathbf G}\!\left(z,\mathbf z_E\right)$ is the partition function  defined by \eqref{multiGF}.
\end{lemma}
To see this, expand $\Psi(z,\mathbf z_E)$  into  a sum of monomials in the $m+1$ variables. It follows from \eqref{EulerRocks} that the  coefficient of the monomial $z^k \prod_{uv\in E} z_{uv}$ is the number of solutions of \eqref{MISSlackForm}.  Solutions correspond to  independent sets of order $k$, so in bracket notation we have   
\begin{align}
\label{mypicoef}
    A_k=\Big [z^k\textstyle \prod_{uv\in E} z_{uv} \Big ] \Psi_{\mathbf G}\left(z,\mathbf z_E\right).
\end{align}
Substituting into \eqref{GFis} gives  the  cumbersome    expression 
\begin{align}
    \pi_{\alpha}(z) =\sum_{k=1}^\infty\left\{ \Big [z^k \textstyle \prod_{uv\in E} z_{uv} \Big ] \Psi_{\mathbf G}\left(z,\mathbf z_E\right) \right \}  z^k. \nonumber
\end{align}
Note that 
\begin{align}
    \nonumber
    \Big [z^k \textstyle \prod_{uv\in E} z_{uv} \Big ] \Psi_{\mathbf G}(z,\mathbf z_E) &= \left [z^k \right]\left\{\Big[\textstyle \prod_{uv\in E} z_{uv} \Big ]\Psi_{\mathbf G}(z,\mathbf z_E)\right\}, 
\end{align}
and let $\Phi(z,\mathbf z_E)$ denote the term in braces. Since 
\[
\pi_{\alpha}(z) = \sum_{k\ge 0}\big\{ [z^k]\Phi(z,\mathbf z_E)\big\} z^k = \Phi(z,\mathbf z_E),
\]
Eqn.  \eqref{mypie} follows. 
\qed

\bigskip
As an aside, we point out that for modest size graphs the coefficients $A_k$ can be found recursively using symbolic algebra packages (e.g., Mathematica).  Let $\mathbf z_E=\big(z_{E}(1),\ldots,z_{E}(m)\big)\in \mathbb C^{m}$, 
where, as mentioned above, the  order of the variables is  the same as the order of $\mathbf x_E$.  Nesting the coefficient extraction of \eqref{mypie} gives    
\begin{align}
\label{nesting}
\pi_{\alpha}(z)=\big [z_E(m)\big]\Big( \cdots\big[z_E(2)\big]\Big (\big[ z_E(1)\big ] \Psi_{\mathbf G}(z,\mathbf z_E)\Big )\cdots \Big).
\end{align}
This procedure  postpones expanding  factors  until  it is necessary to do so.

\section{Independence polynomials are ESPs of the vertex monomials}
\label{explicitISform}
 
It is noted in \cite[Sec. 3.5]{GJ} that coefficient extraction is, in general,  a difficult task. The level difficulty depends on the problem.  In the independent set problem,  we surmise from the  products in the  partition function \eqref{multiGF}  that  ESPs can be used to extract the coefficient of   $\Psi_{\mathbf G}(z,\mathbf z_E)$  that corresponds to $\pi_{\alpha}(z)$. This section confirms this supposition by showing that  the coefficients $A_k$ are  elementary symmetric polynomials of the vertex monomials \eqref{vertexmonomial}. 

We first define ESPs for a general vector $\mathbf x=(x_1,\ldots,x_k)\in\mathbb C^k$ and review a few  basic properties.    The ESP of $\mathbf x$ of order $\ell$  is  
\begin{equation}
\label{ESP}
    e_\ell(\mathbf x)=\sum_{1\le i_1<\cdots<i_\ell\le k} x_{i_1}\cdots\, x_{i_\ell},\quad \ell\ge 1.
\end{equation}
For convenience,   define $e_0(\mathbf x)=1$ for $\ell=0$, so that the GF of the  sequence $e_0(\mathbf x),e_1(\mathbf x),\ldots,e_k(\mathbf x)$ is, by elementary algebra,
\begin{align}
\label{elementaryPoly}
     \sum_{\ell=0}^k e_\ell(\mathbf x) z^\ell=\prod_{\ell=1}^k (1+z \,x_\ell).
\end{align}
ESPs  are  computed recursively as follows.  Let $e_0^{\mu}=1$ for $\mu=0,1,\ldots,k$ and $e_{\ell'}^0=0$ for $\ell'=1,\ldots,\ell$. These values initialize the zero-th row and column, respectively, of a $(\ell+1)\times (k+1)$ table.  The remaining entries in the table are computed by the double recursion 
\begin{align}
\label{doublerecursion}
    e_{\ell'}^{k'}=e_{\ell'}^{k'-1}+x_{k'} \,e_{\ell'-1}^{k'-1},\quad \ell'=1,\ldots,\ell,\quad k'=1,\ldots,k.
\end{align}
The computational complexity of filling the table is $O(k\ell)$.  The $\ell+1$ entries in the last column, namely, $e_0^k, e_1^k,\ldots,e_\ell^k$, are equal to, respectively,  $e_0(\mathbf x),e_1(\mathbf x),\ldots,e_\ell(\mathbf x)$.    

With this background in hand, we  now state the main result. 
\newtheorem{thm}{Theorem}
\newtheorem{corollary}{Corollary}
\begin{thm}
\label{maybe}
Let $\mathbf G=(V,E)$ be a general undirected graph.  For $v\in V$, let  $E(v)\subset E$ denote the set of  edges incident to  $v$. Let $\mathbf Z_C$ be the vector of vertex monomials   \eqref{vertexmonomial}, i.e.,  
\begin{align}
\label{vertexmonomials}
\mathbf Z_E=\big(Z_{E(1)},Z_{E(2)},\ldots,Z_{E(n)}\big)\in \mathbb C^{n}.
\end{align} 
Then  the graph polynomial  of  the number of independent sets of $\mathbf G$ is 
\begin{equation}
\label{theprize}
     \pi_{\alpha}(z) =  \sum_{k=0}^{n} \left|e_k^{\scriptscriptstyle\mathrm{Nil}}\left(\mathbf Z_E\right)\right| z^k,
\end{equation}  
where $e_k^{\scriptscriptstyle\mathrm{Nil}}\left(\mathbf Z_E\right)$ denotes the ESP $e_k\left(\mathbf Z_E\right)$ evaluated for nilpotent edge indeterminates of degree 2, i.e.,
\begin{equation}
\label{nilpotent}
    z_{uv}^2=0,\quad \textrm{for all }\; uv\in E,
\end{equation}
and where  $\left | e_k^{\scriptscriptstyle\mathrm{Nil}}\left(\mathbf Z_E\right) \right|$ is the number of monomial terms in $e_k^{\scriptscriptstyle\mathrm{Nil}}\left(\mathbf Z_E\right)$.
\end{thm}

To see  this,  expand the first product in the GF  \eqref{multiGF}  using   \eqref{elementaryPoly}. This gives, since $n=\left|V\right|$,  
\begin{align}
\label{thetestbegins}
    \Psi_{\mathbf G}(z,\mathbf z_E)=\left (\sum_{\ell=0}^{n} e_\ell(\mathbf Z_{E}) z^\ell \right )\left (\prod_{uv\in E}\left(1+z_{uv}\right)\right ).
\end{align}
The coefficient form  \eqref{mypie} of $\pi_\alpha(z)$ is, using  \eqref{thetestbegins},   
\begin{align}
\label{onestepforPeople}
    \pi_{\alpha}(z) &=\Big [{\textstyle \prod_{uv\in E}}\, z_{uv} \Big ] \left (\sum_{\ell=0}^{n} e_\ell(\mathbf Z_{E}) z^\ell \right )\!\left (\prod_{uv\in E}\left(1+z_{uv}\right)\right )\nonumber\\
    &=\sum_{\ell=0}^{n} z^\ell \; \Big [{\textstyle \prod_{uv\in E}}\, z_{uv} \Big ]  \!\left \{\!e_\ell(\mathbf Z_{E})  \prod_{uv\in E}\left(1+z_{uv}\!\right)\right \}.  
\end{align}
Using the ESP definition  \eqref{ESP} and substituting  the explicit form \eqref{vertexmonomial} of the vertex monomials gives  
\begin{align}
\label{simmple}
    e_\ell(\mathbf Z_{E})&=\sum_{1\le i_1<\cdots<i_\ell\le n} Z_{E({i_1})}\cdots Z_{E({i_\ell})}\\
    &=\sum_{1\le i_1<\cdots<i_\ell\le n}\;\prod_{j=1}^\ell  \,\prod_{uv\in E(i_j)} z_{uv} .\nonumber
\end{align}
Let   $\mathcal E(i_{1:\ell})=E(i_1)\cup\cdots\cup E(i_\ell)$ denote the multiset union, i.e., the set union with repeated elements  retained.   Denote  the multiplicity of an edge $uv\in \mathcal E(i_{1:\ell})$ by  $\mu(uv|{i_{1:\ell}})$. For edges $uv\in E$, but not in $\mathcal E(i_{1:\ell})$, define $\mu(uv|{i_{1:\ell}})=0$. This notation extends the multiplicity function $\mu(\cdot|i_{1:\ell})$ to all edges in $E$, so we can write 
\begin{align}
\label{adetail}
    e_\ell(\mathbf Z_{E})=\sum_{1\le i_1<\cdots<i_\ell\le n}\;\prod_{uv\in E} z_{uv}^{\mu(uv|{i_{1:\ell}})}. 
\end{align}
Substituting this expression into \eqref{onestepforPeople} and moving the bracket inside the sum over the indices $i_{1:\ell}$   gives 
\begin{align}
\label{gdterms}
    \pi_{\alpha}(z)&=\sum_{\ell=0}^{n} z^\ell \sum_{1\le i_1<\cdots<i_\ell\le n}\;   \Big [{\textstyle \prod_{uv\in E}}\, z_{uv} \Big ] \left ( \prod_{uv\in E} z_{uv}^{\mu(uv|{i_{1:\ell}})} \left(1+z_{uv}\!\right)\right ).
\end{align}
Expanding the product in large parentheses gives a sum of $m+1$ terms, each of which is a monomial in the indeterminates  $\{z_{uv}:uv\in E\}$. Terms that have an indeterminate with multiplicity $\mu(uv|{i_{1:\ell}})\ge 2$ are not in the coefficient of $\prod_{uv\in E} z_{uv}$. The indeterminate variables satisfy the nilpotent condition  \eqref{nilpotent}, so these terms, and only these terms, are equal to zero.  The remaining terms have indeterminates with multiplicities   $\mu(uv|{i_{1:\ell}})\le 1$. Since the indeterminates are nilpotent, we have        
\begin{align}
\label{hardtowrite}
      &\Big [{\textstyle \prod_{uv\in E}}\, z_{uv} \Big ] \left ( \prod_{uv\in E} z_{uv}^{\mu(uv|{i_{1:\ell}})} \left(1+z_{uv}\!\right)\right ) = 1. 
\end{align}
The sum over the indices $i_{1:\ell}$  in \eqref{gdterms} is therefore  equal to the number of terms in $e_\ell(\mathbf Z_{E})$ when it is evaluated for nilpotent indeterminates of degree 2. The number of such terms is, by definition,  $\left| e_\ell^{\scriptscriptstyle\mathrm{Nil}}(\mathbf Z_E)\right|$.  
\qed 

\bigskip 
Note that the upper limit $n$ in \eqref{theprize} can be replaced by $\alpha(\mathbf G)$, the size of the maximum independent set.  
\begin{corollary} \label{cor1}
    The monomials in the symmetric function $e_k^{}(\mathbf Z_E)$, when evaluated under the nilpotent condition \eqref{nilpotent},  correspond one-to-one to independent sets of order $k$. 
\end{corollary}
To see this, replace $\ell$ with $k$ in     \eqref{simmple} to obtain  $e_k(\mathbf Z_E)$ as a sum of products of vertex monomials $Z_{E({i_1})}\cdots Z_{E({i_k})}$. Since the indeterminates are nilpotent, it follows from \eqref{vertexmonomial} that these products are nonzero if and only if 
\[
E(i_1)\cap \cdots \cap E(i_k) = \varnothing.
\]
The intersection is empty if and only if the vertices in the set union $E(i_1)\cup\cdots\cup E(i_k)$ are independent sets. 
\qed
\bigskip

Said less formally, to compute the number of independent sets,  write  $Z_{E(v)}$ as a product of   the edge indeterminates $z_{uv}$ and implement the ESP recursion   \eqref{doublerecursion} for the product 
\[
\textstyle \prod_{v\in V } \left(1+z\, Z_{E(v)}\right).
\]
At the end of each  step  of the recursion,  simplify the result  using the nilpotent rule $z_{uv}^2=0$. The number of terms in the simplified expression of the coefficient of $z^k$ is $A_k$.

\begin{corollary} 
    If the edge indeterminates are  nilpotent  \eqref{nilpotent}, then the set of $n$ vertex monomials that comprise the vector $\mathbf Z_E$ in \eqref{vertexmonomials} form a nilpotent algebra whose nilpotent index,  $\eta(\mathbf G)$, is  
    \begin{align}
    \label{nilpotentindex}
        \eta(\mathbf G)=1+\alpha(\mathbf G),
    \end{align}
    where $\alpha(\mathbf G)$  is the order of a maximum independent set.   
\end{corollary}
An algebra for which there is an integer $k\ge 1$ such that the product of any $k$ elements is zero is, by definition, nilpotent \cite{encyc}.  The vertex monomials is a nilpotent algebra, as is  seen by taking $k$ equal to the order of the graph, $n$.   The  index of nilpotency, $\eta(\mathbf G)$, is the smallest $k$ for which this is true.  The vertices  in a maximum independent set  map to vertex monomials in $\mathbf Z_E$ and, by Corollary \ref{cor1}, their product is nonzero.   It follows that $\alpha(\mathbf G)=\eta(\mathbf G)-1$.    
\qed
\bigskip

\begin{corollary} 
    The problem of computing the nilpotent index of the algebra of vertex monomials $\mathbf z_E$ of a general undirected graph $\mathbf G$ is NP-complete.   
\end{corollary}
\noindent The problem of computing the order $\alpha(\mathbf G)$ of a maximum independent set  is NP-complete \cite{MISnpcomplete}. The corollary follows immediately.  \qed

\section{Cliques}
\label{Scetioncliques}
Let $C_k$ denote the number of cliques of size $k\ge 0$  in the graph $\mathbf G=(V,E)$. The GF of this sequence is the graph polynomial,   
\begin{align}
    \pi_{\gamma}(z)=\sum_{k=0}^\infty C_k z^k, 
\end{align}
is derived in this section. The derivation is straightforward because a vertex set is independent if and only if it is a clique in the graph complement \cite[Chap. 7]{BondyMurty}. It follows immediately that the GF of the number of cliques in  $\mathbf G$ is the GF of the number of independent sets in the graph complement, $\overline{\mathbf G}$.  

The complement $\overline{\mathbf G}$ of the simple graph $\mathbf G=(V,E)$ is the simple graph with the same vertex set $V$ but with the edge set $\overline E=E_K\setminus E$, where $E_K$ is the edge set of the complete graph $K_n$. Two vertices in $\overline{\mathbf G}$ are adjacent if and only if they are not adjacent in $\mathbf G$.  

The  number of edges in $\overline{\mathbf G}$ is 
\[
\overline m=\left|E_K\setminus E\right|=\tfrac{1}{2}n(n-1) - m.
\]
Let  $\mathbf z_{\overline E}\in\mathbb C^{\overline m}$ denote the  edge indeterminates $\{z_{uv},uv\in\overline E\}$, arranged in the some  fixed  order.   Paralleling \eqref{vertexmonomial},  the vertex monomials of $\overline{\mathbf G}$ are  products of edge indeterminates, 
\begin{align}
    Z_{\overline E(v)}=\prod_{uv\in \overline E(v)} z_{uv},\;\text{ for } v\in V,
\end{align}
where $\overline E(v)$ is the set of edges in $\overline E$ that are incident to $v$.  The vector of vertex monomials of $\overline{\mathbf G}$ is   
\begin{align}
\label{GCvertexmonomials}
\mathbf Z_{\overline{ E}}=\big(Z_{\overline E(1)},Z_{\overline E(2)},\ldots,Z_{\overline E(n)}\big)\in \mathbb C^{n}.
\end{align} 
With this notation, the  $\mathcal E(i_{1:\ell})=E(i_1)\cup\cdots\cup E(i_\ell)$  leading to the univariate GF of  the number of independent sets of $\overline{\mathbf G}$  is, using \eqref{multiGF},  
\begin{align}
\label{multiGFclique}
\Psi_{\overline{\mathbf G}}(z,\mathbf z_{\overline E})&=\prod_{v\in V}\left(1+z  \,  Z_{\overline E(v)}\right)\!\prod_{uv\in \overline E}\left(1+z_{uv}\right).
\end{align}

\begin{corollary}
     The graph polynomial  of the number of cliques of $\mathbf G=(V,E)$ is 
     \begin{align}
\label{mypiCrust}
    \pi_{\gamma}(z) =\Big [\textstyle \prod_{uv\in \overline E} \,z_{uv} \Big ] \Psi_{\overline{\mathbf G}}(z,\mathbf z_{\overline E}).
\end{align}
\end{corollary}  
\noindent This follows immediately from \eqref{mypie} in Lemma \ref{lemma2}. 
\qed

\bigskip
From  Theorem \ref{maybe}, it is seen that \eqref{mypiCrust} can be computed via ESPs.

\section{Vertex covers} 
\label{vertexcovers}
A vertex cover for $\mathbf G$  is a subset of  $V$ that includes at least one endpoint of every edge in $E$.  The order of a vertex cover is its cardinality. Let $B_k$ denote the number of covers of order $k\ge 0$. The GF of the sequence $B_0, B_1, \ldots$ is 
\begin{align}
\label{GFvc}
 \pi_\beta(z)=\sum_{k=0}^\infty B_k z^k.
\end{align}
Since $V$ is a cover, the $ \pi_\beta(z)$ is a polynomial of degree at most $n$.  A minimum vertex cover is a cover of smallest  order. Thus, from  \eqref{GFvc},   
\begin{align}
   \beta(\mathbf G)=\min_k \;\{k : B_k >0 \}
\end{align}
is the covering number, i.e., the order of the minimum cover.     
 
The  complement (in $V$) of a vertex cover is an independent set, and conversely \cite{Gallai1959}.  The coefficients of $ \pi_\beta(z)$  are therefore the same as those of   $\pi_\alpha(z)$, but in reverse order. We have the following result. 

\begin{thm} 
\cite{AkbariOboudi} The graph polynomial  of  the number of vertex covers of an undirected  graph $\mathbf G=(V,E)$ is 
\begin{equation}
     \pi_\beta(z) =  z^{n\,} 
 \pi_\alpha\!\left(z^{-1}\right), \label{theorem1}
\end{equation}
where $\pi_{\alpha}(z)$ is the graph polynomial for the number of independent sets. 
\end{thm}  
\noindent To see this, note that  $B_k=A_{n-k}$.  Multiply both sides of this equation by $z^k$ and sum from $k=0$ to $k=n$.  The left hand side of the sum is $ \pi_\beta(z)$, and  the right hand side is  equal to \eqref{theorem1}.   
\qed
\bigskip

The well-known result  $\alpha(\mathbf G)+\beta(\mathbf G)=n$ follows  from \eqref{theorem1}.  Since both $\pi_\alpha(\mathbf G)$ and $\pi_\beta(\mathbf G)$ are polynomials,  the lowest powers of $z$ on both sides must be equal.  Thus, $z^{\alpha(\mathbf G)}=z^n z^{-\beta(\mathbf G)}$. 

The GF $\pi_\beta(\mathbf G)$ can also be derived from the partition function of the BLP system that characterizes the number of vertex covers. This system is identical to \eqref{MISSlackForm} but with \eqref{MISSlackForm2b} replaced by 
\[
x_u+x_v+x_{uv}=2.
\]
This derivation gives a  multivariate GF  that is identical to  \eqref{multiGF}, i.e., the partition functions for vertex covers and independent sets are the same.          
\begin{corollary} 
The graph polynomial  of  the number of independent sets of an undirected  graph $\mathbf G=(V,E)$ is, from \eqref{mypie}, 
\begin{align}
\label{vpretty}
    \pi_\alpha(z)=\Big [\textstyle \prod_{uv\in E} z_{uv} \Big ] \Psi_{\mathbf G}\left(z,\mathbf z_E\right),
\end{align}
and the graph polynomial  of the number of vertex covers is
\begin{align}
\label{vvpretty}
     \pi_\beta(z)=\Big [\textstyle \prod_{uv\in E} z_{uv}^2 \Big ] \Psi_{\mathbf G}\left(z,\mathbf z_E\right),
\end{align}
where  $\Psi_{\mathbf G}\left(z,\mathbf z_E\right)$ is the partition function given by \eqref{multiGF}.
\end{corollary}

\noindent Computing these expressions for the  graph in Fig. \ref{sixnodeexample}   gives
\begin{align}
    \pi_\alpha(z) &= 1 + 6 z + 8 z^2 + 2 z^3\nonumber\\
     \pi_\beta(z)& = 2 z^3+8 z^4+6 z^5+z^6.
\end{align}
Combinatorial interpretations of the coefficients of other monomials of $\Psi_{\mathbf G}(\mathbf z)$, e.g., $\prod_{uv\in E} z_{uv}^3$,  are not explored here.

\section{Bipartite partitions} 
\label{maxcut}
A set of edges in $E$ is a ``cut-set'' if removing them  from the graph partitions it   into two disjoint subgraphs. The size of a cut-set is its cardinality.    %The vertices $V$  into two  sets, $L$ and $R$, such that $L\cup R=V$ and $L\cap R=\varnothing$.  
Such partitions are called  bipartite partitions.  They are set partitions and should not be confused with partition functions. 

Given  the graph $\mathbf G=(V,E)$, let $D_k$ denote the number of cut-sets of size $k\ge 0$.  The GF of the sequence is defined by 
\begin{align}
    \pi_\omega(z)=\sum_{k=0}^\infty D_k z^k. 
\end{align} 
For example, the complete graph with three  vertices, denoted by $K_3$. If the cut-set is empty, i.e., if $m=0$, we obtain the trivial  bipartite partition consisting of the empty graph $\varnothing$ with no edges and  $K_3$.  There are three possible cut-sets of size $m=1$, but none of them partition $K_3$.  There are three possible cut-sets of size $m=2$, and each corresponds to a  bipartite partition.  There is one possible cut-set of size  $m=3$, but it   separates $K_3$ into three parts, not two.  Thus,   the GF  of the number of cut-sets of $K_3$  is $\pi_\omega(z)=1+3 z^2$.

Let $(L,R)$ be a bipartite partition of $V$, where $L\cap R=\varnothing$ and $L\cup R=V$. Let $x_u,\,u\in V$, be a binary variable that is equal to 0 if  $u\in L$  and equal to 1 if $u\in R$. Let  $\oplus$  denote the exclusive-or (XOR), known colloquially as  the ``one or the other but not both'' operator, so that  
\begin{align}
x_u\oplus x_v\equiv \left\{
\begin{array}{rl}
0 & \quad\textrm{ if }  x_u=0,\,  x_v=0\\
1 & \quad\textrm{ if }  x_u=0,\,  x_v=1\\
1 & \quad\textrm{ if }  x_u=1,\,  x_v=0\\
0 & \quad\textrm{ if }  x_u=1,\,  x_v=1.
\end{array}\right.
\end{align}   
The binary  variable $x_{uv},\,uv\in E$, is equal to 1 of the edge $uv$ is in the cut-set, and 0 if not. The  cut-set is $\{uv\in E\mid x_{uv}=1\}$, and its size is  $\sum_{uv\in E} x_{uv}$.  The number of cut-sets of size $k$ is  the number of solutions of the  system: \begin{subequations} 
\label{BipartiteSolutions}
\begin{align}
&\textstyle \sum_{uv\in E} x_{uv} = k\label{BipartiteSolutionsaa}\\
&\qquad \quad x_{uv}= x_u \oplus x_v,\,\, \textrm{for all }  uv\in E, \label{BipartiteSolutionsbb}\\
\textrm{where}\quad &x_v \in \{0,1\}\,\, \textrm{and}\,\,
x_{uv}\in \{0,1\}.\label{BipartiteSolutionscc}
\end{align}
\end{subequations}
\noindent The XOR equations  \eqref{BipartiteSolutionsbb} are nonlinear, but they can be replaced by a system of binary linear equations.  To see  this, first note that  $x_{uv}=x_u\oplus x_v$ is equivalent to the system of inequalities
\begin{align*}
x_{uv}&\le x_u+x_v\\
x_{uv}&\ge x_u - x_v\\
x_{uv}&\ge -x_u + x_v\\
x_{uv}&\le 2 - x_u - x_v.
\end{align*}
To convert them to a system of linear binary equations with the same solution set, we add  binary ``slack'' variables that we write in the form   $2(1-s_{uv})$, $2(1-t_{uv})$, $2(1-y_{uv})$, and $w_{uv}$, where $s_{uv},t_{uv}$, and $y_{uv}$ are also binary variables  in $\{0,1\}$.  With these slack variables, the system is 
\begin{align*}
x_u&+x_v-x_{uv} +2 s_{uv} \quad\quad\quad =2\\
-x_u&+x_v+x_{uv} + \quad 2 t_{uv} \quad\quad=2\\
x_u&-x_v+x_{uv} +\quad\quad 2 y_{uv} \quad=2\\
x_u&+x_v+x_{uv} +\quad\quad\quad 2 w_{uv} =2.
\end{align*}
The only binary solutions $(x_u,x_v,x_{uv},s_{uv},t_{uv},y_{uv},w_{uv})$ of the system are 
\begin{align*}
(0,0,0,1,1,1,1)\nonumber \\
(0,1,1,1,0,1,0)\nonumber \\
(1,0,1,1,1,0,0)\nonumber \\
(1,1,0,0,1,1,0)\nonumber
\end{align*}
as is easily verified.   Thus, the nonlinear XOR system of $m+1$  equations \eqref{BipartiteSolutions} in $n+m$ variables is equivalent to a linear system with  $4 m+1$   equations in $n+5 m$ variables: 
\begin{subequations} 
\label{BipartiteSolutionsBinary}
\begin{align}
&\;\textstyle \sum_{uv\in E} x_{uv} = k\label{BipartiteSolutionsBinarya}\\
\nonumber
\textrm{and, for all  $uv\in E$}&,\\x_u&+x_v-x_{uv} +2 s_{uv} \quad\quad\quad =2\label{BipartiteSolutionsBinaryb}\\
\label{BipartiteSolutionsBinaryc}-x_u&+x_v+x_{uv} + \quad 2 t_{uv} \quad\quad=2\\
\label{BipartiteSolutionsBinaryd}x_u&-x_v+x_{uv} +\quad\quad 2 y_{uv} \quad=2\\
\label{BipartiteSolutionsBinarye}x_u&+x_v+x_{uv} +\quad\quad\quad 2 w_{uv} =2\\
\textrm{where }\; &x_v,
x_{uv},s_{uv},t_{uv},y_{uv},w_{uv}\in\{0,1\}.\nonumber
\end{align}
\end{subequations}
Define the variables $\mathbf x_V=\{x_{v},v\in V\}$ and $\mathbf x_E=\{x_{uv},uv\in E\}$,  each arranged in some arbitrary but fixed order. Similarly define the slack variables $\mathbf s_E$, $\mathbf t_E$, $\mathbf y_E$, and $\mathbf w_E$.   Let    
\begin{align}
\label{xVxExSetc}
    \mathbf x=(\mathbf x_V,\mathbf x_E,\mathbf s_E,\mathbf t_E,\mathbf y_E,\mathbf w_E). 
\end{align}   
The matrix version of \eqref{BipartiteSolutionsBinarya}-\eqref{BipartiteSolutionsBinarye}  is $A^{(4m+1)\times (n+5m)} \tilde{\mathbf x}^{n+5m} = \mathbf b^{4m+1}$.   The vector $\tilde{\mathbf x}$ is the transpose of the row vector \eqref{xVxExSetc}, and the vector  
$\mathbf b=(k,2,\ldots,2)^T$. The system matrix is   
\begin{align}
\label{edgeadjencymatrixMaxCut}
A&=\left[ 
\begin{array}{c|c|c|c|c|c}
    \textrm - & \textbf 1^{1\times m} & \textrm - & \textrm - & \textrm - &\textrm{-} \\
    \mathbf C^{m\times n} &-\mathbf I_m & 2\mathbf I_m & \textrm - & \textrm - & \textrm - \\
    \mathbf D^{m\times n} & \mathbf I_m &  \textrm - & 2\mathbf I_m & \textrm - & \textrm - \\
    -\mathbf D^{m\times n} & \mathbf I_m &  \textrm - & \textrm - & 2 \mathbf I_m & \textrm - \\
    \mathbf C^{m\times n} & \mathbf I_m &  \textrm - & \textrm - & \textrm - & 2 \mathbf I_m \\
    \end{array}  
    \right]\!,
\end{align}
where  ``$\,\textrm -\,$'' indicates zero entries, $\textbf 1^{1\times m}$ is the row vector with $m$ ones, and $\mathbf I_m$ is the $m\times m$ identity matrix.  The submatrix $\mathbf C =[c_{ij}]$ is the $m\times n$ edge-to-vertex incidence matrix that is defined by $c_{ij}=1$ if vertex $j$ is an endpoint of edge $i$, and $c_{ij}=0$ if it is not. The sum of row $i$ is always 2, since an edge has exactly 2 endpoints.  The sum of column $j$ is the degree of vertex $j$. The submatrix $\mathbf D=[d_{ij}]$ is identical to $\mathbf C$ except that in every row the first $1$ is changed to $-1$. More carefully, in row $i$, let $j_1$ and $j_2$ be the columns for which $c_{ij_1}=c_{ij_2}=1$.  Then $d_{ij_1}=-1$ if and only if $j_1<j_2$.  The sign change is due to the linear ordering of the vertices.  

We assign the indeterminate   $z\in \mathbb C$ to \eqref{BipartiteSolutionsBinarya}. We need $m$ indeterminates for each   equation in \eqref{BipartiteSolutionsBinaryb}-\eqref{BipartiteSolutionsBinarye};  they are denoted    
\begin{align}
    \mathbf z_{\scriptscriptstyle \textrm {XOR}}= \left(\mathbf z_{1,\scriptscriptstyle E},\,\mathbf z_{2,\scriptscriptstyle E},\,\mathbf z_{3,\scriptscriptstyle E},\,\mathbf z_{4,\scriptscriptstyle E}\right)\in\mathbb C^{4m},
\end{align}
where   $\mathbf z_{j,\scriptscriptstyle E}=\{ z_{j,uv}:uv\in E\}\in\mathbb C^m$ for $j=1,2,3,4$. Paralleling the derivation for independent sets \eqref{EulerRocks}, we obtain
\begin{align}
    %\label{wearenotEuler}
    &\prod_{uv\in E} \Big\{z^{x_{uv}} \; z_{1,uv}^{x_u+x_v-x_{uv} +2 s_{uv}}\;
     z_{2,uv}^{-x_u+x_v+x_{uv} + 2 t_{uv}}\nonumber\\
    &\quad\qquad \times z_{3,uv}^{x_u-x_v+x_{uv} + 2 y_{uv}} \;z_{4,uv}^{x_u+x_v+x_{uv} +2 w_{uv} }\Big\}=z^k\prod_{uv\in E} \prod_{j=1}^4 z_{j,uv}^2 .\nonumber
\end{align}
The partition function is defined as the sum of the left hand side of this equation over the variables in $\mathbf x$, which are listed in \eqref{xVxExSetc}.  Rearranging terms in the products, as was done in the independent sets problem, is tedious but straightforward for all the variables in $\mathbf x$ except $\mathbf x_V$. Calculating gives 
\begin{align}
\label{GrandCut}
    \Psi_{\mathbf G}(z,\mathbf z_{\scriptscriptstyle \textrm {XOR}})&=\Psi_E(z,\mathbf z_{\scriptscriptstyle \textrm {XOR}})\, \Psi_{\scriptscriptstyle \textrm {Slack}}(\mathbf z_{\scriptscriptstyle \textrm {XOR}}) \,\Psi_V(\mathbf  z_{\scriptscriptstyle \textrm {XOR}}),
\end{align}
where 
\begin{align} 
\label{partE}
\Psi_E(z,\mathbf z_{\scriptscriptstyle \textrm {XOR}})&=\prod_{uv\in E}  \left(1+z \,z_{1,uv}^{-1} \, z_{2,uv}\, z_{3,uv} \,z_{4,uv}\right)\\
\Psi_{\scriptscriptstyle \textrm {Slack}}(\mathbf z_{\scriptscriptstyle \textrm {XOR}})&=\prod_{uv\in E}\prod_{j=1}^4\left(1+z_{j,uv}^2\right)\\
\label{partV}\Psi_V(\mathbf z_{\scriptscriptstyle \textrm {XOR}})&= \sum_{\mathbf x_V} \,z_{1,uv}^{x_u+x_v}\,z_{2,uv}^{-x_u+x_v}\,z_{3,uv}^{x_u-x_v}\,z_{4,uv}^{x_u+x_v}.
\end{align}
The final but key step is to write the summand of  $\Psi_V(\mathbf z_{\scriptscriptstyle \textrm {XOR}})$ as a separable product of  vertex monomials. Because of the minus signs, the vertex monomials  are  Laurent monomials, i.e., they are products of indeterminates with powers in the integers $\mathbb Z$. (Laurent polynomials also arise naturally in the combinatorics of knot theory \cite[p. 388]{GodsilRoyle}.)   Define the vector $[\,\ln \mathbf z_{\scriptscriptstyle \textrm{XOR}}]\in\mathbb C^{4m\times 1}$ to be the component-wise natural logarithm of $\mathbf z_{\scriptscriptstyle \textrm{XOR}}$.  Define  $Q\in\mathbb C^{1\times n}$ by  
\begin{align}
\label{Q}
Q=\left[\,\ln \mathbf z_{\scriptscriptstyle \textrm{XOR}}\right]^T \left[ \begin{array}{c}
    \mathbf C\\
    \mathbf D  \\
    \!\!-\mathbf D \\
    \mathbf C  \\
    \end{array} 
    \right].
\end{align}
Let $Q=\big[q(v)\big]_{v\in V}$. The Laurent vertex monomials are 
\begin{align}
\label{Laurentvertexmonomial}
    Z_{E(v)}=\mathrm e^{q(v)}\;\text{ for } v\in V.
\end{align}
Each monomial $Z_{E(v)}$ is a product of integer powers of the $4m$ indeterminate variables  $\{ z_{j,uv}: j\in\{1,2,3,4\},\;uv\in E\}$.  Calculation gives 
\begin{align}
\Psi_V(\mathbf z_{\scriptscriptstyle \textrm {XOR}})&= \prod_{v\in V} (1+Z_{E(v}).
\end{align}
On substituting this expression into \eqref{GrandCut}, it is seen that  $\Psi_{\mathbf G}(z,\mathbf z_{\scriptscriptstyle \textrm {XOR}})$ is a product of Laurent monomials in $4m+1$  indeterminates. 

\begin{thm} The graph polynomial of the number of bipartite cuts of the undirected graph $\mathbf G=(V,E)$ is the coefficient 
\begin{align}
\label{cutAnswer}
    \pi_\omega(z) = \left[\prod_{uv\in E} z_{1,uv}^{2}\,z_{2,uv}^{2}\,z_{3,uv}^{2}\,z_{4,uv}^{2}\right] \Psi_{\mathbf G}(z,\mathbf z_{\scriptscriptstyle \textrm {XOR}}),
\end{align}
where the partition function is given by \eqref{GrandCut}.
\end{thm}
\qed

\bigskip
The Laurent monomial can be replaced by the usual monomials by multiplying through by an appropriate product of the indeterminates. This product is  
\begin{align}
    \phi=\prod_{uv\in E} z_{1,uv}^{}\,z_{2,uv}^{}\,z_{3,uv}^{}.
\end{align}
The variables $z_{4,uv}$ are absent from $\phi$, as one might expect since they appear with  negative exponents in \eqref{partE} and \eqref{partV}. To compensate for multiplying the partition function by $\phi$, we extract a different coefficient:
\begin{align}
\label{cutAnswer2}
    \pi_\omega(z) = \left[\prod_{uv\in E} z_{1,uv}^{3}\,z_{2,uv}^{3}\,z_{3,uv}^{3}\,z_{4,uv}^{2}\right] \Big(\phi \, \Psi_{\mathbf G}(z,\mathbf z_{\scriptscriptstyle \textrm {XOR}})\Big).
\end{align}
Further details are omitted.

For the  example depicted in Fig. \ref{sixnodeexample}, we have 
\begin{align}
D^{7\times 6}=\left[\begin{array}{cccccc}
-1&1&0&0&0&0\\
-1&0&0&0&0&1\\
0&-1&1&0&0&0\\
0&-1&0&0&0&1\\
0&0&-1&1&0&0\\
0&0&-1&0&1&0\\
0&0&0&0&-1&1
\end{array}
\right]\!.
\end{align}
The matrix $C$ is identical to $D$ except that the $-1$s are replaced by $+1$s.   It is of little value  to list all the Laurent  monomials here, but for those who wish to check the calculations, we give  the first and the last, using the edge ordering of  Sect. \ref{Example1}:
\begin{align}
    Z_{E(12)}&= 
    z_{1,12}^{} z_{1,16}^{} \;
    z_{2,12}^{-1} z_{2,16}^{-1} \;
    z_{3,12}^{} z_{3,16}^{} \;
    z_{4,12}^{} z_{4,16}^{} \nonumber\\
    Z_{E(56)}&=
    z_{1,16}^{} z_{1,26}^{} z_{1,56}^{} \;
    z_{2,16}^{} z_{2,26}^{} z_{2,56}^{} \;
    z_{3,16}^{-1} z_{3,26}^{-1} z_{3,56}^{-1} \;
    z_{4,16}^{} z_{4,26}^{} z_{4,56}^{} ,\nonumber
\end{align} 
where  small spaces are inserted  to ease reading. Using \eqref{cutAnswer2} and Mathematica to compute $\pi_\omega(z)$ via the nesting method  \eqref{nesting}, we obtain
\begin{align}
\label{exnumbomega}
\pi_\omega(z)=1+z+4 z^2+10 z^3 + 9 z^4 + 5 z^5 + 2 z^6.
\end{align}
This expression can be verified by inspection.  

As an aside, it is known \cite[Sect. 3.6]{RV} that if the $n$ vertices of a simple undirected graph are partitioned uniformly at random over all possible $2^n$ partitions, then the expected size of the cut-set is $m/2$, that is, \begin{align*}
    \mathbb E\left[ \,|\textrm{random cut-set}|\,\right]= \frac{m}{2} = \frac{\pi'_\omega(1)}{\pi_\omega (1)},
\end{align*}
where the prime denotes differentiation with respect to $z$.  Calculating the quantity on the right hand side for \eqref{exnumbomega} gives  $\frac{7}{2}$, as anticipated.

\section{Concluding remarks}
\label{Concludingremarks}

This paper derives expressions for graph polynomials for  the numbers of  independent sets, cliques, vertex covers, and bipartite cuts  as coefficients of   multivariate partition functions. The sets, cliques, cover, and cuts  can be enumerated explicitly by incorporating into the partition functions  a new set of appropriately defined indeterminate variables called  labels.    %They   ``come along for the ride'' when using symbolic software such as Mathematica. 
For example, without labelling,   $[z^k]\pi_\alpha(z)$ is the number of independent sets of size $k$; however, the same calculation using a labelled partition function  yields a sum of  ``label monomials.''  These  monomials correspond 1-to-1 to the  independent sets of size $k$.    

The methods of this paper can be applied to weighted problems, provided  the weights are integers.   The  partition function remains a multivariate polynomial, or a Laurent polynomial if the weights are negative.   The  degree of the partition function  increases with the size of the weights, and that of course causes practical difficulties.   The famous Ising problem of physics can be studied by the methods of this paper.  This includes Ising problems with an externally imposed magnetic field.  

Saddle point methods are used in analytic combinatorics in several variables (ACSV) to study asymptotic behavior  \cite{PW}  when the  partition function has   one or more saddle points.  Unfortunately, the partition functions  used in this paper are polynomials and do not satisfy the usual  saddle point conditions (e.g.,  \cite[p. 549]{FS}).  By  extending the sums --- which  we   restrict  to $\{0,1\}$ ---  to infinity, we obtain  partition functions  that  \textit{do} have saddle points. This leads to    Hardy-Littlewood circle style methods \cite{HL} for deriving asymptotic estimates of the number of solutions of Diophantine systems.    For example, the partition function  \eqref{neededattheend} of the example becomes a product of geometric series, namely    
\begin{equation}
\label{attheend}\nonumber
    \begin{split}
    \Psi(z,\mathbf z_E)&=(1-z\, z_{12} z_{16})^{-1} (1-z\, z_{12} z_{23} z_{26})^{-1}\cdots(1-z\, z_{16} z_{26} z_{56})^{-1} \\
    &\qquad \times (1-z_{12})^{-1}(1-z_{16})^{-1}\cdots(1-z_{56})^{-1}.
    \end{split}
\end{equation} 
This  function  has saddle points. In the language of ACSV, it has amoebas galore \cite[Sec. 7.3]{PW}. 

\bigskip 
\noindent  \textit{Acknowledgements.}  The author thanks Dr. James Ferry (Metron) for verifying Theorem 1 in Mathematica, and Dr. Kevin Bongiovanni (Raytheon) and Dr. Thomas Wettergren (University of Rhode Island) for helpful discussions.

\end{document}